\documentclass[12pt]{article}
\usepackage{amsfonts}

\newcommand{\oper}[2]{\newcommand{#1}{\mathop{\mathrm{#2}}\nolimits}
}

\oper{\tr}{tr}
\oper{\adj}{adj}
\oper{\Div}{div}
\oper{\ad}{ad}
\oper{\Ad}{Ad}
\oper{\End}{End}
\oper{\Hom}{Hom}
\oper{\Aut}{Aut}
\oper{\SO}{SO}
\oper{\SP}{Sp}
\oper{\SU}{SU}
\oper{\GL}{GL}
\oper{\T}{T}
\oper{\U}{U}
\oper{\id}{I}
\oper{\ext}{Ext}
\oper{\rank}{rank}
\oper{\diag}{Diag}

\def\bR{{\bf R}}

\def\ci{{\cal I}}

\def\cl{{\cal L}}

\def\cu{{\cal U}}
\def\cv{{\cal V}}

\def\homega{\hat{\omega}}

\newcommand{\lie}[1]{\mathfrak{#1}}

\newtheorem{proposition}{Proposition}
\newtheorem{theorem}{Theorem}
\newtheorem{definition}{Definition}

\newcommand{\bproof}{\noindent{\it Proof: }}
\newcommand{\eproof}{\  q.~e.~d. \vspace{0.2in}}

\begin{document}

\title{Reduction of HKT-Structures}

\author{
Gueo Grantcharov
\thanks{
Permanent Address: Department of Mathematics, Sofia
University, 5
"James Bourchier", 1126 Sofia, Bulgaria; Present
address:
Department of Mathematics, University of Connecticut,
Storrs, CT
06269, USA. E-mail: geogran@math.uconn.edu . } \and
George Papadopoulos
\thanks{Address: Department of Mathematics, King's
College London,
Strand, London WC2R 2LS, UK. E-mail:
gpapas@mth.kcl.ac.uk .} \and
Yat Sun Poon
\thanks{
Address: Department of Mathematics, University of
California at
Riverside, Riverside, CA 92521, U.S.A.. E-mail:
ypoon@math.ucr.edu. } }

\maketitle

\noindent{{\bf Abstract}: KT-geometry is the geometry
of a
Hermitian connection whose torsion is a 3-form.
HKT-geometry is
the geometry of a hyper-Hermitian connection whose
torsion is a 3-form. We identify non-trivial
conditions for a
reduction theory for these types of geometry.}

\section{Introduction}

Symplectic reduction is a novel method of constructing
symplectic
 manifolds from others that admit a group action of
symplectic 
diffeomorphisms.
 To describe the main result, let $G$ be a compact
group of symplectic 
diffeomorphisms
 acting on the symplectic manifold $(M,\omega)$ and
$\mathfrak{g}$ be 
the Lie algebra
 of $G$. It can be shown that under certain conditions

$N=\nu^{-1}(\zeta)/G$
 is also a symplectic manifold, where $\zeta\in
\mathfrak{g}^*$ and
 $\nu:~M\rightarrow \mathfrak{g}^*$ is the moment map.
The manifold $N$
 is also denoted with $M/\! /G$.  It is remarkable
that symplectic 
reduction can be
 generalized in various ways. First, it can be shown
that if $M$ is a 
K\"ahler
 manifold admitting a G-action of holomorphic
isometries, then $M/\! 
/G$ is also
 a K\"ahler manifold. Furthermore, it can be shown
that if $M$ is a 
hyper-K\"ahler
 manifold admitting an G-action of tri-holomorphic
isometries, then
 $M/\!/G=\nu^{-1}(\zeta)$ is also hyper-K\"ahler where
 $\nu=(\nu_1, \nu_2,\nu_3): M \rightarrow
\bR^{3}\otimes 
\mathfrak{g}^*$
 and $\zeta=(\zeta_1, \zeta_2,\zeta_3)\in
\bR^{3}\otimes 
\mathfrak{g}^*$ \cite{HKLR}.
 In the context of hyper-K\"ahler reduction there are
three moment maps 
each
 associated to the three complex structures.
One common feature of all symplectic, K\"ahler and
hyper-K\"ahler
reductions is that moment maps exist because the
G-action
preserves some symplectic form.

 More generally it has been shown that if $M$ is a
hypercomplex 
manifold admitting
a tri-holomorphic group action, then $M/\!/G$ is also
hypercomplex  \cite{Joyce1}.
 The details of this construction will be summarized 
in Section 
\ref{HKT-map}. Here
 it is worth mentioning that in the context of
hypercomplex reductions, 
moment
 maps do not arise naturally because in the generic
case there are no
 symplectic forms which are preserved by the group
action. Instead it 
is
 assumed that one can find such functions on $M$ which
have the 
required
 properties.

In the next section, we assume the existence of a
$G$-moment map
on $M$ and study the geometry on the reduced space
$N$. The aim
is to prove that the reduction of a KT-space is a
KT-space and
the reduction of a HKT-space is again a HKT-space. The
definition
and twistor construction of HKT spaces have been given
in
\cite{howepap}. The  properties of KT and HKT
manifolds have been
widely investigated in the literature
\cite{howepap,gibstegp, GP}.
The result on KT-space in Section 2 is not surprising
because a 
Hermitian
structure can easily be found on a reduced space and
every
Hermitian structure has a unique KT-connection. The
existence of
HKT-connection on the reduction of a HKT-space is less
trivial.
Examples of HKT-reduction in this regard are given at
the end of this 
paper.

In the third section, we identify non-trivial and
sufficient
topological or cohomological constraints on either the
manifold
$M$ or the group $G$ to ensure the existence of a
$G$-moment map
on strong KT-manifolds and strong HKT-manifolds. In
the absence
of symplectic forms, this is a non-trivial result as
one usually
generates moment map through the K\"ahler form. In the
fourth
section, we discuss when a potential function on a
HKT-space may
descend to a potential function on the reduced
HKT-space.

\section{Existence of HKT-Structures on Reduced
Spaces}

Assuming the existence of ``moment maps", we examine
the geometry on 
the reduced space in the next two sections.

\subsection{KT Reduction}\label{KT-map}
Before we explain HKT reduction, it is instructive to
consider
first the reduction of KT manifolds, i.e. Hermitian
manifolds
equipped with the hermitian connection whose torsion
is a
three-from.

Let $M$ be a KT manifold and  let $G$ be a compact
group of
complex isometries on $M$. Denote the algebra of
holomorphic
vector fields by $\mathfrak{g}$. Next introduce a
$G$-equivariant
map  $\nu:M\longrightarrow \mathfrak{g}$ satisfying
the
transversality condition, ie $I d\nu(X)\neq 0$ for all
$X\in
\mathfrak{g}$. We remark that a map $\nu$ is
equivariant if
$\nu(g\cdot x)= Adg^* \big(\nu(x)\big)$.

\begin{definition}
 A map $\nu$ is called $G$-moment map if and only if
(i) it is 
equivariant and (ii) it satisfies
the transversality condition.
\end{definition}

 We remark that for simply connected  K\"ahler
manifolds the moment map can be constructed using  the
invariance
of K\"ahler form and complex structure and it
satisfies the
transversality property. However additional conditions
are
required in order the moment map to be equivariant.

 Next given a point $\zeta\in \mathfrak{g}$,
denote the level set $\nu ^{-1}(\zeta )$ by $P$. Since
the map
$\nu $ is $G$-equivariant, level sets are invariant if
the group
$G$ is Abelian or if the point $\varsigma $ is
invariant. Assuming
that the level set $P$ is invariant, and the action of
$G$ on $P$
is free, then the quotient space $N=P/G$ is a smooth
manifold.
Let $\pi :P\to N$ be the quotient map.

It can be shown that in fact $N=P/G$ is a complex
manifold. This
construction can be done as follows. For each point
$m$ in the
space $P$, its tangent space is
\[
T_{m}P=\{t\in T_{m}M:d\nu (t)=0\}.
\]
Consider the vector subspace
\[
\cu_{m}=\{t\in T_{m}P:Id\nu(t)=0\}.
\]
Due to the transversality condition, this space is
transversal to
the vectors generated by elements in $\mathfrak{g}$. 
In
addition, this space is a vector subspace of $T_{m}P$
with
co-dimension $\dim \mathfrak{g}$, and hence it is a
vector
subspace of $T_{m}M$ with co-dimension $2\dim
\mathfrak{g}$. The
same condition implies that, as a subbundle of
$TM_{|P}$, $\cu$
is closed under $I$. Moreover there is a G-invariant
splitting
\begin{equation}
TP = \cu \oplus \cv
\end{equation}
where $\cv$ is the tangent space to the orbits of G
and it is the
bundle of kernels of $d\pi$. We use the terms
``horizontal" and
``vertical" for $\cu$ and $\cv$.

As the projection $\pi$ is an isomorphism on $\cu$,
for any
tangent vector $\hat{A}$ at $\pi (m)$, there exists a
unique
element $A^u$ in ${\cu}_{m}$ such that $d\pi
(A^u)=\hat{A}$. We call
$A^u$ horizontal lift of $\hat{A}$. The complex
structure on $N$
is defined by
\begin{equation}\label{complex}
I\hat{A}=d\pi (IA^u),\quad \mbox{ i.e. } \quad (I{\hat
A})^u
=IA^u.
\end{equation}

\begin{theorem}\label{KT-reduction}
Let $(M,\ci,g)$ be a KT-manifold. Suppose that G is a
compact
group of complex isometries admitting a G-moment map
$\nu$. Then
the complex reduced space $N=M/\!/G$ inherits a KT
structure.
\end{theorem}
\bproof To show this, it suffices to find a complex
structure $I$
and a hermitian metric $g$ on $N$ which are induced
from $M$
because for every Hermitian structure $(I,g)$, there
always
exists a unique KT structure on $N$ \cite{Gauduchon}
\cite{GP}.

To begin, since $\cu$ is $G$-invariant, if $X^u$ is
tangent to
$P$ at $m$ and is contained in $\cu$, then for any
element $f \in
G$, $dL_f(X^u)$ is tangent to $P$ at $f(m)$ and is
contained in
$\cu$. Using $\pi\circ L_f=\pi$, if $X^u$ is a
horizontal lift of
$\hat X$ to a point $m$, then $d\pi \circ
dL_f(X^u)=d\pi
(X^u)=\hat{X}$. Therefore, $dL_f(X^u)$ is the
horizontal lift of
$\hat{X}$ to $f(m)$.

Since $G$ is also a group of isometries, $g(dL_f(X),
dL_f(Y))=g(X, Y)$ for any vectors $X$ and $Y$ tangent
to $P$.
Define a metric $\hat{g}$ on $N$ by
\begin{equation}
{\hat{g}}_{\pi (p)}(\hat{X}, \hat{Y})=g_p(X^u, Y^u)
\end{equation}
where $X^u$ and $Y^u$ are the horizontal lifts of
$\hat{X}$ and
$\hat{Y}$ respectively. From the analysis above, the
metric $\hat
g$ is independent from the choice of the reference
point $p$ of
the orbit. Note that the \lq\lq horizontal" and \lq\lq
vertical"
spaces ARE NOT necessarily orthogonal.

To prove that $\hat{g}$ is Hermitian, we note that
\begin{eqnarray}
g_{\pi(p)}(I{\hat{X}}, I{\hat{Y}})&=&
g_p((I{\hat{X}})^u,
(I{\hat{Y}})^u) =g_p(I({\hat{X}}^u), I({\hat{Y}})^u)
\nonumber\\
&=&g_p({\hat{X}}^u, {\hat{Y}}^u) =g_{\pi (p)}(\hat{X},
\hat{Y}).
\end{eqnarray}
\eproof

\subsection{HKT Reduction}\label{HKT-map}

We shall begin with a description of hypercomplex
reduction
developed by Joyce \cite{Joyce1}. Let $G$ be a compact
group of
hypercomplex automorphism on $M$. Denote the algebra
of
hyper-holomorphic vector fields by $\mathfrak{g}$.
Suppose that
$\nu =(\nu _{1},\nu _{2},\nu _{3}):M\longrightarrow
\bR^{3}\otimes \mathfrak{g}$ is a $G$-equivariant map
satisfying
the following two conditions. The Cauchy-Riemann
condition:
$I_{1}d\nu _{1}=I_{2}d\nu _{2}=I_{3}d\nu _{3}$, and
the
transversality condition: $I_a d\nu_a(X)\neq 0$ for
all $X\in
\mathfrak{g}$. In analogy with a similar definition
given in the
previous section, any map satisfying these conditions
is called a
$G$-moment map. Given a point $\zeta =(\zeta
_{1},\zeta_{2},\zeta
_{3})$ in $\bR^{3}\otimes \mathfrak{g}$, denote the
level set
$\nu ^{-1}(\zeta )$ by $P$. Assuming that the
level set $P$ is invariant, and the action of $G$ on
$P$ is free,
then the quotient space $N=P/G$ is a smooth manifold.

Joyce proved that the quotient space $N=P/G$ inherits
a natural
hypercomplex structure \cite{Joyce1}. His construction
runs as
follows. For each point $m$ in the space $P$, its
tangent space is
\[T_{m}P=\{t\in T_{m}M:d\nu _{1}(t)=d\nu _{2}(t)=d\nu
_{3}(t)=0\}.\]
Consider the vector subspace
\begin{equation}\label{u-distribution}
\cu_{m}=\{t\in T_{m}P:I_{1}d\nu
_{1}(t)=I_{2}d\nu_{2}(t)=I_{3}d\nu _{3}(t)=0\}.
\end{equation}
Due to the transversality condition, this space is
transversal to
the vectors generated by elements in $\mathfrak{g}.$
Due to the
Cauchy-Riemann condition, this space is a vector
subspace of
$T_{m}P$ with co-dimension $\dim \mathfrak{g}$, and
hence it is a
vector subspace of $T_{m}M$ with co-dimension $4\dim
\mathfrak{g}$.

The same condition implies that, as a subbundle of
$TM_{|P}$,
$\cu$ is closed under $I_a$. Moreover there is a
G-invariant splitting
\begin{equation}
TP = \cu \oplus \cv
\end{equation}
where $\cv$ is the tangent space to the orbits of G
and it is the
bundle of kernels of $d\pi$. Again, we use the terms
``horizontal" and ``vertical" for $\cu$ and $\cv$
although these
two spaces are not necessarily orthogonal. Following
techniques and 
notations
of the last section, 
a hypercomplex structure on $N$
is defined by
\begin{equation}\label{hypercomplex}
I_a\hat{A}=d\pi (I_aA^u),\quad \mbox{ i.e. }\quad
(I_aA)^u
=I_aA^u.
\end{equation}

\begin{theorem}\label{HKT-reduction}
Let $(M,\ci,g)$ be a HKT-manifold. Suppose that G is a
compact
group of hypercomplex isometries admitting a G-moment
map $\nu$.
Then hypercomplex reduced space $N=M/\!/G$ inherits a
HKT
structure.
\end{theorem}
\bproof
Define hypercomplex structures $I_a$ on $N=P/G$ as in
(\ref{hypercomplex}). As in the previous section,
define a metric
$\hat{g}$ on $N$ by
\begin{equation}
 g_p(X^u, Y^u)={\hat{g}}_{\pi (p)}(\hat{X}, \hat{Y})
\end{equation}
where $X^u$ and $Y^u$ are the horizontal lifts of
$\hat{X}$ and
$\hat{Y}$ respectively. This is a hyper-Hermitian
metric.

On $M$, define $F_a({{X}}, {{Y}})=g(I_a{{X}}, {{Y}})$
and
\begin{equation}
\omega_1 = F_2 - iF_3.
\end{equation}
This is a (0,2)-form with respect to $I_1$. Since the
hyper-Hermitian structure on $X$ admits a HKT-metric,
${\overline\partial}\omega_1=0$. Equivalently, the
$(0,3)$-part
of $d\omega_1$ vanishes.

Similarly, we define $\homega_1$ on $N$. By
\cite[Proposition
2]{GP}, the hyper-Hermitian metric $\hat{g}$ is a
HKT-metric if
and only if ${\overline\partial}\homega_1=0$. In other
words, 
we need to prove that the
type $(0,3)$-part of $d\homega_1$ with respect to
$I_1$ vanishes.
This is equivalent to
\begin{equation}
\pi^*d{\homega_1}(X^u,Y^u,Z^u) = 0
\end{equation}
for any vectors $X^u,Y^u,Z^u$ in $\cu^{0,1}_{I_1}$. As
\begin{equation}
\pi^* {\homega_1}(Y^u, Z^u) = \omega_1(Y^u, Z^u)
\end{equation}
and we have the following  computation:
\begin{eqnarray*}
& \ &d\pi^*{\homega_1}(X^u,Y^u,Z^u)\\
&=&X^u(\pi^* {\homega_1}(Y^u,Z^u))- Y^u(\pi^*
{\homega_1}(Z^u,X^u))
+ Z^u(\pi^* {\homega_1}(X^u,Y^u))\\
&&-\pi^*{\homega_1}([X^u,Y^u],Z^u)-\pi^*{\homega_1}([Y^u,Z^u],X^u)-
\pi^*{\homega_1}([Z^u,X^u],Y^u)\\
&=&X^u(\omega_1(Y^u,Z^u))-Y^u(\omega_1(Z^u,X^u))+Z^u(\omega_1(X^u,Y^u))\\
& & -\omega_1 ([X^u,Y^u]^u ,Z^u)-\omega_1 ([Y^u,Z^u]^u
,X^u)-\omega_1
([Z^u,X^u]^u,Y^u)\\
&=&d\omega_1(X^u,Y^u,Z^u)\\
& &+\omega_1 ([X^u,Y^u]^v,Z^u) +\omega_1
([Y^u,Z^u]^v,X^u)+\omega_1([Z^u,X^u]^v,Y^u)\\
&=&\omega_1 ([X^u,Y^u]^v,Z^u) +\omega_1
([Y^u,Z^u]^v,X^u)+\omega_1([Z^u,X^u]^v,Y^u).
\end{eqnarray*}

To complete the proof of this theorem we claim that 
$[X^u,Y^u]^v = 0$. Equivalently, $d_a\nu_a ([X^u,Y^u])
= 0$ for
$a=1,2,3$. Since $X^u$ and $Y^u$ are in the kernel of
$d_a\nu_a$
for $a=1,2,3$,
\[
dd_a\nu_a(X^u,Y^u) = X^u(d_a\nu_a(Y^u)) -
Y_a(d_a\nu_a(X^u)) 
-d_a\nu_a([X^u,Y^u]) 
= - d_a\nu_a([X^u,Y^u]).
\]
As $ dd_1\nu_1 $ is of type-(1,1) with respect to
$I_1$ and $X^u
$ and $Y^u$ are type-(0,1) with respect to $I_1$,
$dd_1\nu_1(X^u,Y^u) = 0$. By the Cauchy-Riemann
condition
$d_1\nu_1=d_2\nu_2=d_3\nu_3$, our claim follows.
\eproof

\section{Moment Maps for Strong KT and  HKT-Spaces}
As we have seen, the construction of new HKT manifolds
using HKT
reduction requires the existence of a G-moment map
satisfying the
requirements of Theorem \ref{HKT-reduction}. This
moment map is
not specified within the theory, as it is the case for
the
hyper-K\"ahler reduction, but rather its existence is
an
additional assumption of the construction. However as
we shall
see in the special case of reduction for strong KT (
and  HKT)
manifolds, under certain assumptions, there is such a
moment map
which arises naturally. The local construction of a
moment map
for KT and HKT geometries presented below parallels
the
construction of an action for two-dimensional (2,0)-
and
(4,0)-supersymmetric gauged sigma models with
Wess-Zumino term in
\cite{hullpapspence}, respectively. Again, we focus on
a
reduction theory for strong KT-structure first. The
reduction
theory for strong HKT-structures follow.

\subsection{Local Consideration}
Let $G$ a compact group of complex automorphisms on a
{\it
strong} KT manifold $M$. In particular $G$ is a group
of
isometries on $M$ which leaves in addition the torsion
three-form
$H$ invariant. To continue we introduce a basis
$\{e_a;
a=1,\dots, \dim\mathfrak{g}\}$ in the Lie algebra of
$\mathfrak{g}$ and denote the associated vector fields
of $M$ with
$\{X^a; a=1,\dots, \dim\mathfrak{g}\}$; denote with
$\{e^a;
a=1,\dots, \dim\mathfrak{g}^*\}$ the associated basis
in the dual
$\mathfrak{g}^*$ of $\mathfrak{g}$. The conditions for
invariance
of the KT structure can now be written as
\begin{equation}
{\cal L}_ag=0,~~{\cal L}_aH=0,~~{\cal L}_aI=0
\end {equation}
where ${\cal L}_a={\cal L}_{X^a}$; similarly later for
the inner
derivation we have $i_a=i_{X^a}$.

Using the assumption that $M$ is a strong KT manifold,
$dH=0$,
the last equation above implies that
$di_aH=0$
and so there is a locally defined one-form $u_a$ such
that $$
i_aH=du_a\ . $$ Clearly $u_a$ is uniquely defined up
to the
addition of a closed one-form.

Next let us denote with $\tilde X$ the one-form dual
with the
vector field $X$ with respect to the KT metric. Using 
${\cal
L}_a I=0$, one can show that the two-form
$d(\tilde X_a+u_a)$
is type-(1,1) with respect to the complex structure
$I$. Therefore, by
the ${\overline\partial}$-Poincare Lemma, there is a 
locally defined complex-valued function $h_a$ on $M$
such that
$(\tilde X_a+u_a)^{1,0}=\partial h_a$. Let $f_a$ be
the real part
of $h_a$. Define
\begin{equation}
w_a={\tilde X}_a+u_a-df_a.
\end{equation}
Then $w_a^{1,0}=i\partial \nu_a$ where $\nu_a$ is a
constant
multiple of the imaginary part of $h_a$. Therefore, we
can write
\begin{equation}
w_a= I d\nu_a\ .
\end{equation}
Let $\xi=\xi^ae_a$ be any element in $\lie{g}$. Define
a map $\nu$
from $M$ to $\lie{g}^*$ by
\begin{equation}
\nu(x)(\xi):=\sum_a \xi^a\nu_a(x).
\end{equation}

A necessary condition for  $\nu$ to be well-defined on
$M$ is
that the class of $i_a H$ in $H^2(M,\bR)$ should be
trivial. If
in addition $M$ satisfies the
$\partial\bar\partial$-lemma, then
$\nu$ will be well-defined on $M$.

In the case when the group $G$ is Abelian, the issue
of equivariance is
absent and hence the map $\nu$ so constructed is the
moment map.
Before we investigate equivariance in general, we
consider the issue of 
non-degeneracy.
\begin{definition}
A holomorphic Killing vector field $X$ is
non-degenerate if
$d\nu_X\neq 0$.
\end{definition}
Therefore, a holomorphic Killing vector field is
non-degenerate if
its moment map is non-constant. The following
proposition is
useful to determine when a holomorphic Killing vector
field is
non-degenerate.

\begin{proposition} If the length of a holomorphic
Killing vector field is non-constant, the vector field
is 
non-degenerate.
\end{proposition}
\bproof
Note that $d\nu_X=0$ if and only if $d{\tilde
X}+du=0$.
i.e. $d{\tilde X}+\iota_XH=0$. It means that for any
vector field
$Y$ and $Z$,
\begin{eqnarray*}
Y(g(X,Z))-Z(g(X, Y))-g(X,[Y,Z])+H(X, Y, Z) &=&0\\
\mbox{ i.e. } \hspace{.5in} g(\nabla_YX, Z)+g(X, 
\nabla_YZ)-g(\nabla_ZX, Y) & & \\ 
-g(X, \nabla_ZY)-g(X, [Y,Z])+H((X, Y, Z) &=&0\\
\mbox{ or } \hspace{.5in} \ g(\nabla_YX,
Z)-g(\nabla_ZX, Y)+2H(X, Y, 
Z)&=&0.
\end{eqnarray*}
On the other hand, since $\cl_Xg=0$ and $\nabla g=0$,
\begin{eqnarray*}
0&=& X(g(Y, Z))-g([X,Y], Z)-g(Y,[X,Z])\\
&=&
g(\nabla_XY,Z)+g(Y,\nabla_XZ)-g([X,Y],Z)-g(Y,[X,Z])\\
&=& g(\nabla_YX,Z)+g([X,Y],Z)+H(X, Y,Z)+g(Y,
\nabla_ZX)\\
& &+g(Y, [X,Z])+H(Y,X,Z)-g([X,Y],Z)-g(Y,[X,Z])\\
&=&g(\nabla_YX, Z)+g(Y,\nabla_ZX).
\end{eqnarray*}
Combining the above two identities, we find that for
any vector
fields $Y,Z$,
\begin{equation}
g(\nabla_YX,Z)=-H(X, Y,Z)=-du_X(Y,Z).
\end{equation}
In particular, $g(\nabla_YX,X)=0$ for any $Y$. Since
$\nabla
g=0$. It implies that $dg(X,X)=0$. \eproof

\subsection{Equivariance}
Now we seek conditions for $\nu$ to be equivariant.
This issue
will be analyzed in the next few paragraphs. The map
$\nu$ is 
equivariant if and only if $\nu (g\cdot x)=Ad
g^{*}(\nu (x))$. Let $X$ be any element in $\lie{g}$.
The
equivariance is determined by
\begin{equation}
\nu (g\cdot x) (X) = \nu(x)(Ad g (X)).
\end{equation}
The infinitesimal version of the above identity is
\begin{equation}
{\cal L}_Y\nu_X=\nu_{[V,X]}; \quad \mbox{ equivalently
} \quad
{\cal L}_Y\nu_X-\nu_{[Y,X]}=0.
\end{equation}

Let $[X_b, X_a]=f_{ba}^cX_c$ be the structural
equations for the
algebra $\lie{g}$ so that $f_{ba}^c$ are constants.
Apply the above formula to $w_a$
and $u_a$ respectively with respect to $X_b$, the
equivariance
conditions for $w_a$ and $u_a$ are
\begin{equation}\label{equivariant}
\cl_bw_a-f_{ba}^cw_c=0, \quad {\cal
L}_bu_a-f_{ba}^cu_c=0.
\end{equation}
These are  non-trivial conditions. Note that
\begin{eqnarray}
d{\cal L}_bu_a &=& {\cal L}_bdu_a={\cal L}_b\iota_aH =
\iota_{{\cal L}_bX_a}H+\iota_a{\cal L}_bH=\iota_{{\cal
L}_bX_a}H
\nonumber\\
&=& f_{ba}^c\iota_cH=f_{ba}^cdu_c=d(f_{ba}^cu_c).
\label{closed}
\end{eqnarray}
By Poincar\'e lemma, there exists  a locally defined
closed
1-form $v_{ba}$ such that
\begin{equation}
{\cal L}_bu_a-f_{ba}^cu_c=v_{ba}.
\end{equation}
Therefore, $v_{ba}$ is the obstruction for $u_a$ to be
equivariant.

Next, note that $\cl_ag=0$,
\begin{eqnarray*}
(\cl_b{\tilde X}_a)X &=& \cl_b(g(X_a, X))-g(X_a,
\cl_bX)\\
&=& g(\cl_bX_a, X)+g(X_a, \cl_bX)-g(X_a,\cl_bX)
=f_{ba}^cg(X_c, X)\\
&=& f_{ba}^c {\tilde X}_cX.
\end{eqnarray*}
Therefore, the $\mathfrak{g}^*$-valued 1-form
$w:=w_a\eta^a$ is
equivariant if and only if $u:=u_a\eta^a$ is
equivariant.

Assuming that $u$ is equivariant. This implies that
after a
possible shift of $u_a$ with respect to a closed
one-form, $u_a$
must satisfy the above equation. Note that even if
 $u_a$ is equivariant, it is not unique but rather
defined up to an 
equivariant
{\it closed} one-form.

Next since $dw_a$ is an (1,1)-form and if we assume
that the
$\partial\bar\partial$-lemma  applies on the manifold
$M$ (see
either \cite[5.11]{DGMS} or \cite[Corollary
2.110]{Besse}), there
is a function $\nu_a$ on M such that
\begin{equation}\label{dw}
dw_a=dd^cv_a=dId\nu_a.
\end{equation}
Therefore, the 1-form
\[
z_a=w_a-Id\nu_a
\]
is closed. In the above equation $\nu_a$ is not
uniquely defined
but rather it is defined up the addition of the real
part of a
{\it holomorphic } function.

As we have assumed that $u_a$ is equivariant, $w_a$ is
equivariant.  We 
obtain
\begin{equation}\label{thonca}
Id\nu_{ba}+z_{ba}=0
\end{equation}
where $\nu_{ba}={\cal L}_b\nu_a-f_{ba}{}^c \nu_c$ and
$z_{ba}={\cal L}_bz_a-f_{ba}{}^c z_c$. Since
$dz_{ba}=0$,
(\ref{thonca}) implies that $dd^c\nu_{ba}=0$. By
$\partial{\overline\partial}$-Lemma again, $\nu_{ba}$
is a
harmonic function and hence is the real part of a
holomorphic
function $f_{ba}$. If, in addition,
\begin{equation}\label{dF}
f_{ba}={\cal L}_b F_a-f_{ba}{}^c F_c
\end{equation}
for some holomorphic functions $F_a$, then redefining
$\nu_a$ as
$\nu_a-{\rm Re}F_a$ and $z_a$ as $z_a-d{\rm Im}F_a$
both $\nu_a$
and $z_a$ become equivariant. So  there is a choice of
$u_a$,
such that $w_a=Id\nu_a$. Therefore, we have found an
equivariant
moment map $\nu: M\rightarrow \mathfrak{g}^*$.

\subsection{Cohomology}
The various conditions that we have found for the
existence of a
moment map in the previous section can be identified
as classes
in de-Rham $H^*_{dR}$ and in $H^*_\delta$ cohomology,
where
$\delta$ will be defined shortly. Let $\delta_G$ be
the map
defining Lie algebra cohomology in the usual way
\cite{GS}. In
particular, for $\theta\in \mathfrak{g}^*$ and $\zeta,
\eta\in
\mathfrak{g}$,
\begin{equation}
\delta_G\theta(\zeta,\eta)=-\theta([\zeta,\eta]).
\end{equation}
Therefore, (note the convention for wedge product) in
terms of
structural constants with respect to the dual basis
$\theta^a$,
\begin{eqnarray}
\delta_G\theta^a &=&
-\sum_{b,c}f_{bc}^a\theta^b\otimes\theta^c =
-\sum_{b<c}f_{bc}^a(\theta^b\otimes\theta^c-\theta^c\otimes\theta^b)
\nonumber\\
&=& -\sum_{b<c}f_{bc}^a\theta^b\wedge\theta^c =
-\frac12
\sum_{b,c}f_{bc}^a\theta^b\wedge\theta^c.
\label{convention}
\end{eqnarray}
In particular, $\delta^2_G=0$. Next for $\phi$ in
$\Lambda^\ell(M)$ and $X$ in $\mathfrak{g}$,  define
\begin{equation}
{\hat\delta}\phi(X):=\cl_X\phi.
\end{equation}
Equivalently, ${\hat\delta}\phi= {\cal L}_a
\phi\cdot\theta^a.$
Then we extend this operator to
\begin{equation} \delta:~~
\Lambda^\ell(M)\otimes
\Lambda^k\mathfrak{g}^*\otimes\mathfrak{g}^*\rightarrow
\Lambda^\ell(M)\otimes\Lambda^{k+1}\mathfrak{g}^*\otimes\mathfrak{g}^*
\end{equation}
 as follows. If $\phi$ is in $\Lambda^\ell(M)$,
$\theta$ is in
$\Lambda^{k}\mathfrak{g}^*$ and $\eta$ is in
$\mathfrak{g}^*$,
then define
\begin{eqnarray}
\delta(\phi\cdot \theta\otimes \eta) :=
{\hat\delta}\phi\wedge\theta\otimes
\eta+\phi\cdot\delta_G\theta\otimes\eta
+(-1)^k\phi\cdot\theta\wedge\delta_G\eta.
\end{eqnarray}
This map generates a resolution.
\[
\Lambda^\ell(M)\otimes\mathfrak{g}^*\stackrel{\delta_0}{\rightarrow}
\Lambda^\ell(M)\otimes
\Lambda^1\mathfrak{g}^*\otimes\mathfrak{g}^* \cdots
{\rightarrow}
\Lambda^\ell(M)\otimes
\Lambda^k\mathfrak{g}^*\otimes\mathfrak{g}^*
\stackrel{\delta_k}{\rightarrow}
\Lambda^\ell(M)\otimes
\Lambda^{k+1}\mathfrak{g}^*\otimes\mathfrak{g}^*
\cdots.
\]
We claim that this resolution is a complex. i.e.
$\delta_k\circ\delta_{k+1}=\delta^2=0$. To check,
notice that
\begin{eqnarray*}
& & \delta^2(\phi\cdot \theta\otimes\eta)\\
&=& \delta({\hat\delta}\phi\wedge\theta)\otimes \eta
+(-1)^{k+1}{\hat\delta}\phi\wedge\theta\wedge\delta_G\eta
+{\hat\delta}\phi\wedge\delta_G\theta\otimes\eta
+\phi\cdot \delta^2_G\phi\otimes\eta\\
& &+(-1)^{k+1}\phi\cdot \delta_G\theta\wedge\delta_G\eta
+(-1)^k{\hat\delta}\phi\wedge\theta\wedge\delta_G\eta
+(-1)^k\phi\cdot\delta_G\theta\wedge\delta_G\eta\\
& &
+(-1)^{2k}\phi\cdot\theta\wedge\delta_G^2\eta\\
&=& \delta({\hat\delta}\phi\wedge\theta)\otimes \eta+
{\hat\delta}\phi\wedge\delta_G\theta\otimes\eta =
\left(\delta({\hat\delta}\phi\wedge\theta)+
{\hat\delta}\phi\wedge\delta_G\theta \right)\otimes\eta\\
&=&
\left(\delta({\hat\delta}\phi)\wedge\theta-{\hat\delta\phi}\wedge\delta_G\theta+
{\hat\delta}\phi\wedge\delta_G\theta\right)\otimes\eta
= \left(\delta(\cl_a\phi \theta^a)\right)\wedge\theta\otimes\eta\\
&=&
\left(\cl_b\cl_a\phi\cdot\theta^b\wedge\theta^a
+\cl_c\phi\cdot\delta_G\theta^c\right)
\wedge\theta\otimes\eta
= \left(\cl_b\cl_a\phi-\frac12 f^c_{ba}\cl_c\phi\right)\cdot
\theta^b\wedge\theta^a\wedge\theta\otimes\eta.
\end{eqnarray*}
Since $[\cl_a,\cl_a]\phi=f^c_{ba}\cl_c\phi$ and
$f^c_{ba}=-f^c_{ab}$,
\begin{equation}
\cl_b\cl_a\phi-\frac12f^c_{ba}\cl_c\phi
=\cl_a\cl_b\phi+\frac12f^c_{ba}\cl_c\phi
=\cl_a\cl_b\phi-\frac12f^c_{ab}\cl_c\phi.
\end{equation}
It shows that the term $\cl_b\cl_a\phi-\frac12
f^c_{ba}\cl_c\phi$
is symmetric in the indices $ab$ while the term
$\theta^b\wedge\theta^a$ is skew symmetric in $ab$. It
follows
that $\delta(\cl_a\phi \theta^a)=0$ and hence
$\delta^2=0$ as
claimed.

One can now define a cohomology theory with respect to
$\delta$
in the usual way and denote it with
\begin{equation}
H_\delta^k(\Lambda^\ell(M)\otimes\mathfrak{g}^*)
:=\frac{\ker\delta_k}{\mbox{\rm image} ~
\delta_{k-1}}.
\end{equation}
Since $\delta$ commutes with $d$, one can also
naturally define
the cohomology groups
$H_\delta^k(C^\ell(M)\otimes\mathfrak{g}^*)$, where
$C^\ell(M)$
are the close $\ell$-forms on $M$.

A cohomology theory based on a resolution of  ${\cal
O}\otimes
\mathfrak{g}^*$, where ${\cal O}$ is the sheaf of
germs of
holomorphic functions on $M$, is similarly defined.
This is
possible because the group $G$ consists of holomorphic
actions.
In particular,
${\overline\partial}\circ\cl_a=\cl_a\circ{\overline\partial}$.
This cohomology is
\begin{equation}
H_\delta^k({\cal O}\otimes\mathfrak{g}^*)
:=\frac{\ker\delta
:{\cal
O}\otimes\Lambda^k\mathfrak{g}^*\otimes\mathfrak{g}^*
}
{\mbox{\rm image}~ \delta: {\cal
O}\otimes\Lambda^{k-1}\mathfrak{g}^*\otimes\mathfrak{g}^*}.
\end{equation}

Returning now in the discussion of the previous
section, we
have seen a necessary condition for the existence of a
moment map
in the KT case is that $i_aH$ is a trivial class in
$H^2_{dR}(M)$. Now we write $i_aH=du_a$ and   define
$u=u_a \eta^a.$
This is a section of $\Lambda^1(M)\otimes
\mathfrak{g}^*$. Using
(\ref{convention}),
\begin{eqnarray}
\delta u&=& (\delta u_a)\otimes\eta^a+\sum_c
u_c\delta_G\eta^c
=\cl_b u_a\theta^b\otimes\eta^a-\sum_{a,b,c}f^c_{ba} 
u_c\theta^b\otimes\eta^a\\
&=& \sum_{a,b}\left( \cl_b u_a-\sum_c
f^c_{ba}u_c\right)\theta^b\otimes\eta^a.
\end{eqnarray}

Due to (\ref{closed}), the 1-form part is closed.
Therefore
$\delta u$ is an element of
$C^1(M)\otimes\Lambda^1\mathfrak{g}^*\otimes\mathfrak{g}$.
Obviously, it is in the kernel of $\delta$. It defines
a class in
$H_\delta^1(C^1(M)\otimes\mathfrak{g}^*)$. Since $u$
is not
necessarily a closed 1-form, this class is not
necessarily
trivial although it is represented by $\delta u$. Due
to
computation of previous paragraphs, this cohomology
class is the
obstruction for adjusting $u$ by a closed 1-from so
that it could
be equivariant.

If this class vanishes, then as we have explained
$\delta (w_a
\eta^a)=0$ as well. Using this and assuming that
$\nu_a$ is
well-defined in $w=Id\nu+z$, where $\nu=\nu_a\eta^a$
and
$z=z_a\eta^a$, we have $Id\delta\nu+\delta z=0$. As we
have
explained in the previous section the obstruction for
both $z$
and $\nu$ to be equivariant are  $\delta z$ and
$\delta\nu$
respectively. The last identity implies that it
suffices to find
the condition for $\delta \nu=0$.

Due to identity (\ref{dw}), $\delta w=0$ and
$\cl_aI=0$, we have
$dId\delta\nu=0$. Therefore, by
$\partial{\overline\partial}$-Lemma, there exists
holomorphic
function $f_{ba}$ such that $\nu_{ba}={\rm Re}
f_{ba}$. Define
\begin{equation}
f:=f_{ba}\theta^b\otimes\eta^a.
\end{equation}
This is an element in ${\cal
O}\otimes\Lambda^1\mathfrak{g}^*\otimes\mathfrak{g}^*$.
The
function part of $\delta f$ is holomorphic as the
group $G$
consists of holomorphic actions.

However, the real part of $\delta f$ is equal to
$\delta \nu=0$.
Therefore, $\delta f$ is purely imaginary. This is
possible only
if $\delta f=0$. It follows that $f$ defines  a class
in
$H_\delta^1({\cal O}\otimes\mathfrak{g}^*)$. Note that
the class
of $f$ vanishes if and only if the equation $f=\delta
F$ has a 
solution. 
In other words,
there are solutions for  the equation (\ref{dF}).

Some of the conditions that we have derived above can
be cast
into an elegant form using equivariant cohomology
\cite{atyhbott}.
In physics, it is known that the obstructions for
gauging {\it
bosonic} two-dimensional  sigma models with
Wess-Zumino term
\cite{hullspence, jjmo} are elements of equivariant
cohomology
\cite{jose}. The theorem below provides sufficient
conditions for
KT reduction.

\begin{theorem}
Let $M$ be a strong KT manifold and $G$ be a compact
group acting
on $M$ and leaving invariant the KT structure. If the
torsion
three-form $H$ admits an equivariant extension as a
closed form in
$EG\times_G M$, $H_\delta^1( {\cal
O}\otimes\mathfrak{g}^*)=0$ and
the $\partial\bar \partial$-lemma applies on $M$, then
$M/\!/G$
is a KT manifold.
\end{theorem}
\bproof Note that $EG$ is the universal classifying
bundle space for 
the
group $G$.
It can be shown that a closed three-form $H$ in $M$
admits
an equivariant extension in $EG\times_G M$,  if $H$ is
invariant
under the group action of $G$ on $M$ and there are
equivariant
one-forms $\{u_a; a=1,\dots, {\rm dim}\mathfrak{g}\}$
on $M$ such
that
 \begin{equation}\label{equiv}
i_aH=du_a ~~{\rm and}~~ i_au_b+i_bu_a=0\ .
\end{equation}
 Of course the one-form $u_a$ is defined up to the
addition of an equivariant closed one-form $v_a$.
Because of
this, the one-form $w_a=u_a+\tilde X_a$ is equivariant
and $dw_a$
is an (1,1) form on $M$. If the
$\partial\bar\partial$-lemma
applies, then $w_a=Id\nu_a+z_a$, where $\nu_a$ is a
function on $M$
and $z_a$ is closed one-form. It can be shown that in
fact $z_a$
is equivariant. Indeed, since $w_a$ is equivariant and
the
G-action preserves the complex structure, we have
\begin{equation}\label{thonc}
Id\nu_{ba}+z_{ba}=0
\end{equation}
where $\nu_{ba}={\cal L}_b\nu_a-f_{ba}{}^c \nu_c$ and
$z_{ba}={\cal L}_bz_a-f_{ba}{}^c z_c$. We have seen
that the
obstruction for $z_a$ and $\nu_a$ to be equivariant
lies in
$H^1_\delta({\cal O}\otimes \mathfrak{g}^*)$. Since
this vanishes
$z_a$ and $\nu_a$ are equivariant. So there is a
choice of $u_a$,
such that $w_a=Id\nu_a$.

It remains to prove the transversality condition. This
follows
from the last condition in (\ref{equiv}) because it
implies that
$i_a u_b$ is skew-symmetric and so $i_a w_b$ is the
sum of a
non-degenerate symmetric matrix with $i_a u_b$.
Therefore $\nu$
is a G-moment map and so $M/\! /G$ is a KT manifold.
\eproof

\subsection{Moment Maps on Strong HKT Structures}
The construction of  G-moment maps for the reduction
of strong
HKT manifolds can proceed as in the case of strong KT
manifolds
above. The only difference is that for each complex
structure
$\{I_r; r=1,2,3\}$ one gets
\begin{equation}\label{shktmom}
w_a={I_r}d(\nu^r)_a+z^r_a
\end{equation}
where $z^r_a$ are again equivariant closed one-forms
provided
that the obstructions in $H^1_\delta({\cal
O}\otimes\mathfrak{g})$ vanish. In this case however
it is not
always possible to redefine $u_a$ such that
$w_a={I_r}d(\nu^r)_a$
unless $z^1_a=z^2_a=z^3_a$. Nevertheless, we can still
use the map
$\nu: M\rightarrow \bR^3\otimes\mathfrak{g}$ as
defined in
(\ref{shktmom}) as a moment map. This moment map is
equivariant
 but neither transversality nor the Cauchy-Riemann 
conditions
generically hold. Thus we have the following theorem:

\begin{theorem} Let $M$ be a strong HKT manifold and
$G$ be a compact 
group
acting on $M$ and leaving invariant the HKT structure.
If the
torsion three-form $H$ admits an extension as a closed
form in
$EG\times_G M$ such that $w_a={I_r}d(\nu^r)_a$ with
$\nu$
equivariant, then $M/\!/G$ is a HKT manifold.
\end{theorem}
\bproof The proof follows from that of reductions of
strong KT
manifolds and that of reductions of weak HKT
manifolds.
\eproof

\section{Potential Functions}
Recall that if $(M,\ci, g)$ is a HKT manifold with
K\"ahler forms
$\omega_a$, a {\it HKT potential} is a function $\rho$
 such
that $2\omega_1 = dd_1\rho+d_2d_3\rho, 2\omega_2 =
dd_2\rho+d_3d_1\rho, 2\omega_3=dd_3\rho+d_1d_2\rho$.
In
this section, we follow the methods in \cite{KS} to
find a
potential function on reduced space. We continue to
use the
notations established in Section \ref{HKT-map}.

\begin{theorem}
Let $(M, \ci, g)$ be a HKT manifold with HKT potential
function
$\rho$. Suppose that G is a compact group of
hypercomplex
isometries leaving $\rho$ invariant with moment map
$\nu =
(\nu_1,\nu_2,\nu_3)$ such that the tangent vectors to
the orbits
of $G$ in $\nu^{-1}(0)$ are in the $\ker (d_a\rho)$,
for $a=1,2,3$. 
Then the
function $\rho$ induces a HKT potential function on
the reduced
space $N=M/\!/G$.
\end{theorem}
\bproof
Let $P:=\nu^{-1}(0)$ and $i: P \rightarrow M$ be the
inclusion
map. Now we first check that
$i^*dd_a\rho|_{\cu}=dd_ai^*\rho|_{\cu}$ where $\cu$ is
defined in
(\ref{u-distribution}). To this end notice that
\[
i^*d\rho(X^u)=d\rho(X^u), \quad
i^*d\rho([X^u,Y^u])=d\rho([X^u,Y^u])
\]
and $i^*I_ad\rho(X^u) =
I_ad\rho(di(X^u))=-d\rho(I_aX^u)$
because $I_ad\rho(X)=-d\rho(I_aX)$. By direct
computations after
restricting on points of $P$  we have:
\begin{eqnarray*}
& &(i^*dd_a\rho)(X^u,Y^u) =di^*d_a\rho(X^u,Y^u)\\
&=&X^u((i^*I_ad\rho)(Y^u)) -Y^u((i^*I_ad\rho)(X^u)) -
i^*I_ad\rho([X^u,Y^u])
\\
&=&-X^u(d\rho(I_aY^u))+Y^u(d\rho(I_aX^u))+d\rho(I_a[X^u,Y^u])
\\
&=&
-X^u(d\rho(I_aY^u))+Y^u(d\rho(I_aX^u))+d\rho(I_a[X^u,Y^u]^u).
\end{eqnarray*}
The last equality is due to $d\rho(I_a[X^u,Y^u]^v) 
=-d_a\rho([X^u,Y^u]^v)=0$. This is true
because $[X^u,Y^u]^v$ is tangent to an orbit of $G$
and the
condition in the theorem. We shall use the same
argument repeatedly and 
implicitly in subsequent computation.

As the map $\rho$ is $G$-invariant, for $x$ in $P$, we
may define
\begin{equation}
\rho_N(\pi(x)):=\rho(x)
\end{equation}
where $\pi$ is the quotient map from $P$ onto $N=P/G$.
In other
words, $\pi^*\rho_N=\rho$. It follows that
\begin{eqnarray*}
& &(\pi^*dd_a\rho_N)(X^u,Y^u) =
d\pi^*d_a\rho_N(X^u,Y^u) \\
&=&
X^u(d_a\rho_N(d\pi(Y^u)))-Y^u(d_a\rho_N(d\pi(X^u))-d_a\rho_N(d\pi([X^u,Y^u]))
\\
&=& -X^u(d\rho_N(I_ad\pi Y^u))+Y^u(d\rho_N(I_ad\pi
X^u))
+d\rho_N(I_ad\pi
[X^u,Y^u])) \\
&=& -X^u(d\rho(I_aY^u))+Y^u(d\rho(I_aX^u)) 
+d\rho_N(d\pi(I_a[X^u,Y^u]^u)) \\
&=& -X^u(d\rho(I_aY^u))+Y^u(d\rho(I_aX^u))
+d\rho(I_a[X^u,Y^u]^u).
\end{eqnarray*}
It follows that
$i^*dd_a\rho_{|\cu}=\pi^*dd_a{\rho_N}_{|\cu}$.
Similarly,
\begin{eqnarray*}
& &(\iota^*d_ad_b\rho)(X^u, Y^u)=
(I_adI_cd\rho)(di(X^u), di(Y^u))\\
&=& (I_adI_cd\rho)(X^u, Y^u)=dI_cd\rho(I_aX^u,
I_aY^u)\\
&=&
I_aX^u(d_c\rho(I_aY^u))-I_aY^u(d_c\rho(I_aX^u))-d_c\rho
([I_aX^u,
I_aY^u])\\
&=&
I_aX^u(d_c\rho(I_aY^u))-I_aY^u(d_c\rho(I_aX^u))-d_c\rho
([I_aX^u, I_aY^u]^u)..
\end{eqnarray*}
On the other hand,
\begin{eqnarray*}
& &(\pi^*d_ad_b\rho_N)(X^u, Y^u)=
(I_adI_cd\rho_N)(d\pi(X^u), 
d\pi(Y^u))\\
&=& (dI_cd\rho_N)(I_ad\pi(X^u), I_ad\pi(Y^u))
=(dI_cd\rho_N)(d\pi I_a(X^u), d\pi I_a(Y^u))\\
&=& (\pi^*dI_cd\rho_N)(I_aX^u, I_aY^u)=
(d\pi^*I_cd\rho_N)(I_aX^u,
I_aY^u)\\
&=&
I_aX^u(\pi^*I_cd\rho_N(I_aY^u)-I_aY^u(\pi^*I_cd\rho_N(I_aX^u))
-\pi^*I_cd\rho_N([I_aX^u, I_aY^u])\\
&=& I_aX^u(-d\rho_N(d\pi
(I_cI_aY^u))-I_aY^u(-d\rho_N(d\pi
(I_cI_aX^u))
-d\rho_N(d\pi I_c[I_aX^u, I_aY^u])\\
&=&
I_aX^u(I_cd\rho(I_aY^u))-I_aY^u(I_cd\rho(I_aX^u))+d\rho(d\pi
I_c[I_aX^u,
I_aY^u]^u)\\
&=& I_aX(d_c\rho (I_aY^u))-I_aY^u(d_c\rho
(I_aX^u))-d_c\rho([I_aX^u,I_aY^u]^u)
\end{eqnarray*}
Therefore,
$i^*d_bd_c\rho_{|\cu}=\pi^*d_bd_c{\rho_N}_{|\cu}$ for
all even permutation $(abc)$ of $(123)$. At the end we
use the
fact that the reduced K\"ahler forms
$\overline{\omega}_a$ are
characterized by the condition
$(\pi^*\overline{\omega}_a)_{|\cu}=
(i^*\omega_a)_{|\cu}$  and conclude that
$2\overline{\omega}_a =
dd_a\rho_N + d_bd_c\rho_N$ \eproof

\noindent{\bf Remark:} In the case of when the torsion
vanishes
the condition in the above theorem is equivalent to
the one
proposed by Kobak and Swann \cite{KS}. In both cases
the crucial
point is to ensure $i^*I_ad\rho=I_ai^*d\rho$. In both
cases
$d\rho(X^v)=0$ since $\rho$ is invariant.

\section{Examples}
It is known that $SU(3)$ admits invariant hypercomplex
structure,
constructed by Joyce. Moreover Pedersen and Poon
\cite{PP2}
considered the deformation of this structure and
succeeded to
represent any "small" deformation as a hypercomplex
reduced space
of the space $S^1 \times S^{11}$ under an appropriate
$S^1$
action. As it is shown in \cite{GP} and \cite{OP}, the
space $S^1
\times S^{11}$ is HKT  and one can check that the
$S^1$-actions
considered in\cite[Section 6.3]{PP2} are
HKT-isometries.
Now according to the theorem of section 2.2 we have:
\begin{theorem} Any small deformation of the invariant
hypercomplex
structure on $SU(3)$ admits a HKT structure.
\end{theorem}

In the rest of this section, we will construct new
HKT-metrics through 
a reduction process. 
We begin with a well-known metric, namely the Taub-NUT
metric.

\subsection{Taub-NUT metric}
We use the notation of \cite{GRG}. Let $ {\cal M} =
\mathbb{H}
\times \mathbb{H} $ with quaternionic coordinates
$(q,w)$. We
identify points $(t,x,y,z)\in R^4$ with a quaternion
$q\in H$:
$q=t+ix+jy+kz.$ 
The (quaternion) conjugate is
$\overline{q}=t-ix-jy-kz$. The flat
metric on $M$ is
\begin{equation}
ds^2_{\mbox{flat}}=dqd{\overline{q}}+dwd{\overline{w}}.
\end{equation}
Using left multiplication of the unit quaternions $i,
j$ and $k$,
we find the hypercomplex structure $I, J$ and $K$ such
that
\begin{eqnarray}\label{hcx structure}
Idt &=& dx, Idx=-dt, Idy=dz, Idz=-dy,\nonumber\\
Jdt&=& dy, Jdx=-dz, Jdy=-dt, Jdz=dx, \nonumber\\
Kdt &=&dz, Kdx=dy, Kdy=-dx, Kdz=-dt.
\end{eqnarray}
With respect to these complex structures, the K\"ahler
form of
the flat metric $dqd{\overline q}$ are
\begin{equation}
\omega_I = dt\wedge dx+dy\wedge dz,
\quad
\omega_J = dt\wedge dy+dz\wedge dx,
\omega_K = dt\wedge dz+dx\wedge dy.
\end{equation}

Let $G$ be $\mathbb{R}$, $t\in \mathbb{R}$ with the
action
$(q,w)\rightarrow (qe^{it},w+\lambda t)$, for
$\lambda$ in
$\mathbb{R}.$
This is a group of hyper-K\"ahler isometries. It
generates a
moment map:
\begin{equation}
\nu = \frac{1}{2}qi\overline{q} +
\frac{\lambda}{2}(w-\overline{w})
\end{equation}
We write ${\bf r} = qi\overline{q}$, $r = |\bf{r}|$
and $w=y +
{\bf y}$ so ${\bf r}$ and $\bf y$ are in $R^3$.
Moreover,
\begin{equation}
\nu = \frac{1}{2}{\bf r} + \lambda {\bf y}.
\end{equation}
Define $\psi$ by $q=\rho e^{i\psi /2}$ where $\rho$ is
a pure
quaternion. Now using the coordinates $(\psi, {\bf r},
y, {\bf
y})$, we write the flat metric on M as
\begin{equation}\label{flat metric}
ds^2_{\mbox{flat}} = \frac{1}{4}(\frac{1}{r}d{\bf r}^2
+ r(d\psi +
{\bf\omega}\cdot d{\bf r})^2) + dy^2 + d{\bf y}^2,
\end{equation}
where ${\mbox{curl}{\bf \omega}}=\mbox{grad}
(\frac{1}{r})$.
In these coordinates the $G$-action is
\begin{equation}\label{action}
(\psi,y) \rightarrow (\psi + 2t, y+ \lambda t),
\end{equation}
which leaves $\tau = \psi - 2y/\lambda$ invariant.  On
$\nu^{-1}(0)$, one has ${\bf
y}=-\frac{1}{2\lambda}{\bf r}$. The
induced metric in the coordinates $({\bf r}, \tau, y)$
on
$\nu^{-1}(0)$ is
\begin{equation}
ds^2_{\mbox{flat}}= \frac{1}{4} (\frac{1}{r}d{\bf r}^2
+ r(d\tau
+\frac{2}{\lambda}dy + {\bf\omega}\cdot d{\bf r})^2) +
dy^2 +
\frac{1}{4\lambda^2}d{\bf r}^2.
\end{equation}

The quotient space $\nu^{-1}(0)/G$ is obtained by an
orthogonal
projection along the Killing vector field
$\partial/\partial y$.
It turns out that the quotient metric  is the Taub-Nut
metric:
\begin{equation}
ds^2_{\mbox{TN}} =
\frac14\left(\frac1{r}+\frac1{\lambda^2}\right)d{\bf
r}^2
+\frac14\left(\frac1{r}+\frac1{\lambda^2}\right)^{-1}
\left(
d\tau+{\bf \omega}\cdot d{\bf r} \right)^2.
\end{equation}

\subsection{A HKT-version of Taub-NUT
metric}\label{hkt reduction}
Given the preparation of the last section, we are now
ready to consider 
HKT-reduction.
Let $h$ be a function of $r$. We consider the metric
on
$\mathbb{H}\backslash\{0\}\times \mathbb{H}$ given by
\begin{equation}
ds^2_h = \frac{h(r)}{q{\overline{q}}}dqd{\overline{q}}
+
dwd{\overline{w}} =\frac{h(r)}{r}dqd{\overline{q}} +
dwd{\overline{w}}.
\end{equation}
As $\frac{h(r)}{r}dqd{\overline{q}}$ is a HKT-metric
on
$\mathbb{H}\backslash\{0\}$ and product of HKT metrics
is again a
HKT metric, $ds^2$ is a HKT metric. Since the
hypercomplex
structure does not change, the group $G$ remains
hypercomplex. It
is again a group of isometries. Therefore, we again
use the
moment maps $\nu$ generated by the action $G$ with
respect to the
hyper-K\"ahler metric $ds^2_{\mbox{flat}}$. On
$\nu^{-1}(0)$ the
induced metric with respect to $ds^2_h$ is
\begin{eqnarray}\label{taubhkt}
& & \frac{h}{4r} (\frac{1}{r}d{\bf r}^2 + r(d\tau
+\frac{2}{\lambda}dy + {\bf\omega}\cdot d{\bf r})^2) +
dy^2 +
\frac{1}{4\lambda^2}d{\bf r}^2
\\
&=& \frac14\left( \frac{h}{r^2}+\frac{1}{\lambda^2}
\right)d{\bf
r}^2 +\left( 1+\frac{h}{\lambda^2} \right)dy^2
+\frac{h}{2\lambda} dy\odot (d\tau+\omega\cdot d{\bf
r})
+\frac{h}4(d\tau+\omega\cdot d{\bf r})^2.\nonumber
\end{eqnarray}
Here we used $\alpha\odot \beta = \alpha\otimes\beta +
\beta\otimes\alpha$. So $\alpha\odot\alpha = 2
\alpha\otimes\alpha$.

As hyper-K\"ahler reduction is also obtained by
orthogonal
projection, the horizontal distribution $\cu$ is
defined by $\ker
\theta$. Therefore, the reduced metric is obtained by
taking the
restriction of $ds^2$ on $\nu^{-1}(0)$ modulo $\theta$
or $\mu$
where
\begin{equation}
\theta=\iota_{\frac{\partial}{\partial
y}}ds^2_{\mbox{flat}},
\quad
\mu=dy+\frac{1}{2\lambda} \frac{( d\tau+{\bf
\omega}\cdot d{\bf
r})}{(\frac1{r}+\frac1{\lambda^2})}.
\end{equation}
In other words, if $\hat{g}$ is the quotient metric,
then there
is a 1-form $\alpha$ and function $a$ on $\nu^{-1}(0)$
such that
\begin{equation}
ds^2_h=a\mu\otimes\mu + (\alpha\otimes\mu
+\mu\otimes\alpha)
+{\hat{g}}.
\end{equation}
It follows that $\iota_{\frac{\partial}{\partial
y}}ds^2=a\mu
+\alpha$. In our example,
\begin{equation}
a= 1+\frac{h}{\lambda^2}, \quad
\alpha = \frac{1}{2\lambda} \left( h-\left(
1+\frac{h}{\lambda^2}\right) \left(
\frac{1}{r}+\frac{1}{\lambda^2}\right)^{-1}\right)
(d\tau+
\omega\cdot d{\bf r}).
\end{equation}
Therefore the quotient metric is
\begin{eqnarray}
& & \frac{1}{4}\left(
\frac{h}{r^2}+\frac{1}{\lambda^2}
\right)d{\bf r}^2 +\left( 1+\frac{h}{\lambda^2}
\right)dy^2
+\frac{h}{2\lambda} dy\odot (d\tau+\omega\cdot d{\bf
r})
+\frac{h}{4}(d\tau+\omega\cdot d{\bf r})^2
\nonumber\\
& &-\left( 1+\frac{h}{\lambda^2}\right)
\left(dy+\frac{1}{2\lambda} \frac{(d\tau+{\bf
\omega}\cdot d{\bf
r})}{\frac{1}{r}+\frac{1}{\lambda^2}} \right)^2
\nonumber\\
& & -\left(dy+\frac{1}{2\lambda} \frac{(d\tau+{\bf
\omega}\cdot
d{\bf r})}{\frac{1}{r}+\frac{1}{\lambda^2}} \right)
\odot
\frac{1}{2\lambda} \left( h-\left(
1+\frac{h}{\lambda^2}\right)
\left(
\frac{1}{r}+\frac{1}{\lambda^2}\right)^{-1}\right)
     (d\tau+ \omega\cdot d{\bf r})
\nonumber\\
&=& \frac{1}{4} \left(
\frac{h}{r^2}+\frac{1}{\lambda^2}
\right)d{\bf r}^2 +\frac{1}{4} \left(
\frac{h}{r^2}+\frac{1}{\lambda^2} \right) \left(
\frac{1}{r}+\frac{1}{\lambda^2} \right)^{-2}
(d\tau+ \omega\cdot d{\bf r})^2\\
&=& \left( \frac{h}{r^2}+\frac{1}{\lambda^2} \right)
\left(
\frac{1}{r}+\frac{1}{\lambda^2} \right)^{-1}
ds^2_{\mbox{TN}}.
\end{eqnarray}
In particular, the quotient metric is conformally
equivalent to the
Taub-NUT metric, a hyper-K\"ahler metric.

Amongst the class of weak HKT metrics that have been
constructed
above, there is a strong HKT metric which is complete.
This is
\begin{equation}\label{strong}
ds^2=\big({1\over r}+{1\over\lambda^2}\big) ds^2_{TN}\
.
\end{equation}
For this metric, the function $h$ is
\begin{equation}
h(r)=1+{2\over \lambda^2} r+{1\over
\lambda^2}\big({1\over
\lambda^2}-1\big) r^2 ~.
\end{equation}
This metric is strong HKT because the conformal factor
is a
harmonic function with respect to the Taub-NUT
hyper-K\"ahler
metric. The asymptotic behavior of the metric is as
follows: As
$r\rightarrow \infty$, the metric (\ref{strong})
approaches the
standard metric on $S^1\times \bR^3$. As $r\rightarrow
0$, the
metric (\ref{strong}) approaches $$ ds^2\sim {1\over
r} \big( r
(d\tau+\omega\cdot d{\bf r})^2+ {1\over r} d{\bf
r}^2\big)\ . $$
Changing back to the quaternionic coordinates $q$, we
find that
the above metric can be rewritten as $$ ds^2={1\over
q\bar q} dq
d\bar q= du^2+ds^2(S^3) $$ with $r=q \bar q$ and
$u=\log(|q|)$.
So it is the standard metric on $\bR\times S^3$. In
physics
language, the metric (\ref{strong}) interpolates
between the
ten-dimensional Kaluza-Klein vacuum and the near
horizon geometry
of the NS5-brane.

\subsection{ A HKT-Version of Lee-Weinberg-Yi Metric}
We are interested in  examples beyond four-real
dimension.
As noted in \cite{GRG}, a high dimension analog of the
Taub-NUT metric
 is the Lee-Weinberg-Yi (LWY) metric. We construct a
family HKT-version 
of this metric.
Moreover, these metrics are not conformal to the
LWY-metric.

We first review the construction of LWY-metric very
briefly to fix 
notations.
We take ${\cal M}={\mathbb{H}}^m\times {\mathbb{H}}^m$
with coordinates $(q_a, w_a)$, $a=1, \dots, m$. Let
$\Lambda=(\lambda_a^b)$ be a real non-degenerate
$m\times m$-matrix. Let $V=(v_a^b)$ be the inverse
matrix.
For $G=\mathbb{R}^m=(t_1, \dots, t_m)$, define
an action by
\begin{equation}
q_a\mapsto q_ae^{it_a}, 
\quad
w_a\mapsto w_a+\sum_b \lambda_a^bt_b.
\end{equation}
With respect to the flat metric
$ds^2_{\mbox{flat}}=\sum_a
dq_ad{\overline{q}}_a+\sum_a 
dw_ad{\overline{w}}_a$
and the hypercomplex structure defined as in (\ref{hcx
structure}), the 
group $G$ is a group hyper-holomorphic
isometries. The moment map
\begin{equation}
\nu=(\nu_1, \dots, \nu_m): {\cal M} \to
{\mathbb{R}}^m\otimes 
{\mathbb{R}}^3
\end{equation}
is given by
\begin{equation}
\nu_a=\frac12 q_ai{\overline{q}}_a +\frac12 \sum_b
   \lambda_a^b(w_b-{\overline{w}}_b).
\end{equation}
Define ${\bf r}_a=q_a i{\overline{q}}_a$, $r_a=|{\bf 
r}_a|=q_a{\overline{q}}_a$,
${\bf {y}}_a=\frac12(w_a-{\overline{w}}_a)$.
It follows that $w_a=y_a+{\bf y}_a$. Now
${\bf{r}}_a$ and ${\bf{y}}_a$ are in $\mathbb{R}^3$
and
the moment map is
\begin{equation}
\nu_a=\frac12{\bf r}_a+\sum_b\lambda_a^b{\bf y}_b.
\end{equation}
Define $\psi_a$ by $q_a=\rho_a e^{i\psi_a/2}$ where
$\rho_a$ is a pure 
quaternion. Now using the coordinates
$(\psi_a, {\bf r}_a, y_a, {\bf y}_a)$, one may
construct explicitly a 
hyper-K\"ahler metric on the
quotient space in the way the Taub-NUT metric is
constructed. This is 
the LWY-metric.

For reference in subsequent computation, we note that
in these 
coordinates the $G$-action is 
$(\psi_a,y_a) \rightarrow (\psi_a + 2t_a, y_a+
\sum_b\lambda_a^b t_b)$.
It leaves the functions
\begin{equation}
\tau_a = \psi_a -2\sum_b v_a^by_b
\end{equation}
invariant.
On the level set $\nu^{-1}(0)$,
${\bf r}_a=-2\sum_b\lambda_a^b{\bf y}_b$.
Equivalently,
${\bf y}_a=-\frac12 \sum_b v_a^b{\bf r}_b.$

Next, consider a new metric on $\cal M$:
\begin{equation}
ds^2 = \sum_a
f_a({q_a{\overline{q}}_a})dq_ad{\overline{q}}_a
+\sum_a dw_ad{\overline{w}}_a
\\
= \sum_a f_a(r_a)dq_ad{\overline{q}}_a
+\sum_a dw_ad{\overline{w}}_a.
\end{equation}
This is a HKT-metric. The group $G$ is again a group
of 
hyper-holomorphic
isometries. We may use the $G$-moment map $\nu$ again
to construct a
quotient metric  $\hat{g}$ with respect to $ds^2$.

The restriction of the metric $ds^2$ on $\nu^{-1}(0)$
with respect to the coordinates $({\bf r}_a, \tau_a,
y_a)$ is
\begin{eqnarray}
 & & \sum_a
\left(
\frac{f_a}{4r_a}d{\bf r}_a^2
+ \frac{f_ar_a}{4}
( 2 \sum_b v_a^b dy_b
+ d\tau_a + A_a)^2
+
 dy_a^2
+\frac{1}{4}(\sum_b v_a^bd{\bf r}_b)^2
\right)
\nonumber
\\
&=&
\frac14
\sum_{b,c}
\left(
\frac{\delta_b^cf_c}{r_c}+\sum_a v_a^bv_a^c
\right) d{\bf r}_b\otimes d{\bf r}_c
+\sum_{b,c}
\left(
\delta_b^c
+
\sum_a (f_a r_a v_a^bv_a^c)
\right) dy_b\otimes dy_c
\nonumber\\
& &
+\frac12 \sum_{a,b} f_ar_av_a^b dy_b\odot(d\tau_a+A_a)
+\frac14\sum_a f_ar_a(d\tau_a+A_a)^2.
\end{eqnarray}

To find the quotient metric $\hat{g}$, it suffices to
find functions 
$F_{ab}$ and
1-forms
$\alpha_a$ such that
\begin{equation}
ds^2=\sum_{a,b}F_{ab}\theta_a\otimes\theta_b
+\sum_a (\theta_a\otimes
\alpha_a+\alpha_a\otimes\theta_a)
+{\hat{g}}.
\end{equation}

Now the problem is that the Killing
vector fields $\frac{\partial}{\partial y_a}$
generated by $G$ on the 
zero
level set in general
are not mutually orthogonal.

{}From now on, we limit our discussion to the case
when 
$\lambda_a^b=\lambda_a\delta_a^b$. 
Equivalently, $\Lambda$  is a diagonal matrix
whose non-zero entry is $\lambda_a$. Its inverse is a
diagonal matrix 
whose
non-zero
entry is $v_a=\frac{1}{\lambda_a}$. In this case,
\begin{equation}
\theta_c:=\iota_{\frac{\partial}{\partial y_c}}ds^2
=(1+r_cv_c^2)dy_c+\frac12r_cv_c(d\tau_c+A_c)
=(1+\frac{r_c}{\lambda_c^2})dy_c+\frac{r_c}{2\lambda_c}(d\tau_c+A_c)
\end{equation}
where $A_c:=\omega({\bf r}_a)\cdot d{\bf r}_a.$ 
Since the vector fields $\frac{\partial}{\partial
y_a}$ are mutually
orthogonal with respect to $ds^2$,
\begin{equation}
\iota_{\frac{\partial}{\partial y_c}}ds^2
=(1+r_cv_c^2)(\sum_a F_{ca}\theta_a+\alpha_c).
\end{equation}
The restriction of the metric $ds^2$ on $\nu^{-1}(0)$
with respect to the coordinates $({\bf r}_a, \tau_a,
y_a)$ is
\begin{eqnarray*}
& & \sum_a
\left(
\frac{f_a}{4r_a}d{\bf r}_a^2
+ \frac{f_a r_a}{4}
( 2v_a dy_a
+ d\tau_a + A_a)^2
+
 dy_a^2
+\frac{v_a^2}{4} d{\bf r}_a^2
\right)
\nonumber
\\
&=&
\sum_a
\left(
\frac{f_a}{4r_a}d{\bf r}_a^2
+ \frac{f_a r_a}{4}
( \frac{2}{\lambda_a} dy_a
+ d\tau_a + A_a)^2
+
 dy_a^2
+\frac{1}{4\lambda_a^2} d{\bf r}_a^2
\right)
\nonumber
\\
&=&
\sum_a(
\frac{1}{4}
(f_ar_a+v_a^2) d{\bf r}^2_a
+(1+f_ar_av_a^2) dy_a^2
+ \frac12 f_ar_av_a dy_a\odot(d\tau_a+A_a)
\\
& &+\frac14 f_ar_a(d\tau_a+A_a)^2
).
\end{eqnarray*}
Therefore,
\begin{eqnarray*}
\iota_{\frac{\partial}{\partial y_a}}ds^2
&=& (1+f_ar_av_a^2)dy_a+\frac12 v_af_ar_a(d\tau_a+A_a)
\\
&=&
\frac{(1+f_ar_av_a^2)}{(1+r_av_a^2)}\theta_a
+\frac12 f_av_ar_a
\left( 1-
\frac{1+f_a r_av_a^2}{(1+r_av_a^2)f_a}\right)
(d\tau_a+A_a)\\
&=&
\frac{(1+f_ar_av_a^2)}{(1+r_av_a^2)}\theta_a
+\frac{v_ar_a}{2(1+r_av_a^2)}( f_a - 1)
(d\tau_a+A_a).
\end{eqnarray*}
It implies that the matrix $(F_{ab})$ is a diagonal
matrix and
\begin{equation}
F_a =F_{aa}=\frac{(1+f_ar_av_a^2)}{(1+r_av_a^2)^2}, 
\quad
\alpha_a =\frac{v_ar_a}{2(1+r_av_a^2)^2}( f_a - 1)
(d\tau_a+A_a).
\end{equation}
Then the quotient metric is
n\begin{eqnarray}
{\hat{g}} &=& ds^2
-\sum_a
\left(
F_a \theta_a\otimes\theta_a +\theta_a\odot\alpha_a
\right)
\nonumber\\
&=&
\frac14
\left(
\frac{f_a}{r_a}+v_a^2
\right) d{\bf r}_a^2
\nonumber\\
& &+
\left(
\frac{f_ar_a}{4} -
\frac{ r_a^2v_a^2(1+f_ar_av_a^2)}{4(1+r_av_a^2)^2}
-\frac{v_a^2r_a^2(f_a-1)}{2(1+r_av_a^2)^2}
\right)
(d\tau_a+A_a)^2
\nonumber\\
&=&
\frac14
\sum_a
\left(\frac{f_a+r_av_a^2}{1+r_av_a^2}\right)
\left(
\left(\frac{1+r_av_a^2}{r_a}\right)d{\bf r}_a^2
+\left(\frac{1+r_av_a^2}{r_a}\right)^{-1}(d\tau_a+A_a)^2
\right).
\end{eqnarray}

When $f_a=1$ for all $a$, we obtain a simple version
of LWY-metric:
\begin{equation}
ds^2_{\mbox{LWY}}
=\frac14
\sum_a
\left(
\left(\frac{1+r_av_a^2}{r_a}\right)d{\bf r}_a^2
+\left(\frac{1+r_av_a^2}{r_a}\right)^{-1}(d\tau_a+A_a)^2
\right).
\end{equation}
This is simple because this metric is a product
metric.

In general, so long as not all the
$\lambda_a=\frac{1}{v_a}$ are equal, 
the quotient metric $\hat{g}$ is 
a HKT-metric. However, it is no longer conformal to
the LWY-metric. 

\section*{Acknowledgments}

 GP is supported by a University Research
Fellowship from the Royal Society. GG is supported by
the university of California at Riverside.
 This work is partially supported
by SPG grant PPA/G/S/1998/00613 and by the European
contract HPRN-CT-2000-00101.


\begin{thebibliography}{01}

\bibitem{Besse} A. Besse. {\em Einstein Manifolds},
Springer-Verlag, 
New York 1987.

\bibitem{DGMS} P. Deligne, P. Griffiths, J. Morgan \&
D. Sullivan. {\em 
Real homotopy
theory of K\"ahler manifolds}, Inventiones Math.{\bf
29} (1975)
245--274.

\bibitem{Gauduchon}  P. Gauduchon. {\em   Hermitian
connections and Dirac operators}, Bollettino U.M.I.,
{\bf   11B}
(1997) 257--288.

\bibitem{GRG} G. Gibbons \& P. Rychenkova \& R. Goto.
{\em 
HyperK\"ahler
quotient construction of BPS monopole moduli spaces},
Commun.
Math. Phys. {\bf 186} (1997) 581--599.

\bibitem{howepap} P.S. Howe \& G. Papadopoulos. {\em
Twistor Spaces for  
HKT manifolds},
Phys. Lett. {\bf B379} (1996) 80-86;  hep-th/9602108.

\bibitem{gibstegp} G.W. Gibbons, G. Papadopoulos \&
K.S. Stelle.
{\em HKT and OKT Geometries on Soliton Black Hole
Moduli Spaces},
Nucl. Phys. {\bf B508}  (1997) 623--658;
hep-th/9706207.

\bibitem{GRG} G. Gibbons \& P. Rychenkova \& R. Goto.
{\em 
HyperK\"ahler
quotient construction of BPS monopole moduli spaces},
Commun. Math.
Phys. {\bf 186} (1997) 581--599.

\bibitem{GP} G. Grantcharov \& Y. S. Poon. {\em
Geometry of 
hyperK\"ahler
connections with torsion}, Commun. Math. Phys. {\bf
213} (2000)
19--37.

\bibitem{GS} V. Guillemin \& S. Sternberg. {\em
Symplectic Techniques 
in Physics},
Cambridge University Press, Cambridge 1990.

\bibitem{HKLR} N.J. Hitchin, A. Karlhede, U.
Lindstr\"om, M. Ro\v cek.
{\em Hyper-K\"ahler metrics and supersymmetry},
Commun. Math.
Phys. {\bf 108} (1987) 535--589.

\bibitem{Joyce1}  D. Joyce. {\em   The hypercomplex
quotient
and quaternionic quotient}, Math. Ann. {\bf   290}
(1991)
323--340.

\bibitem{KS} P. Kobak \& A. Swann. {\em Hyper K\"ahler
potentials via
finite-dimensional quotients}, preprint;
math.DG/0001027.

\bibitem{OP}  A. Opfermann \& G. Papadopoulos.
{\em   Homogeneous HKT and QKT manifolds}, preprint;
math-ph/9807026.

\bibitem{hullpapspence}
 C.M. Hull, G. Papadopoulos \& B. Spence. {\em Gauge
Symmetries
 for (p,q) Supersymmetric Sigma Models}, Nucl. Phys.
{\bf B363} (1991) 
593--621.

\bibitem{atyhbott} M.F.  Atiyah \& R. Bott. {\em The
Moment Map and 
Equivariant Cohomology},
Topology {\bf 23} (1984) 1-28.

\bibitem{hullspence}
 C.M. Hull, B. Spence. {\em The Gauged Non-Linear
Sigma Model with 
Wess-Zumino
  Term},  Phys. Lett. {\bf B232}  (1989) 204.

 \bibitem{jjmo} I. Jack, D.R. Jones, M. Mohammedi \&
H. Osborne. {\em
  Gauging the General Non-Linear Sigma Model with a
Wess-Zumino Term},
  Nucl. Phys. {\bf B332} (1990) 359-379.

\bibitem{jose}
J. M. Figueroa-O'Farrill \& S. Stanciu. {\em
Equivariant
Cohomology and Gauged Bosonic Sigma Models},
hep-th/9407149.

\bibitem{PP2}  H. Pedersen \& Y. S. Poon. {\em
Inhomogeneous hypercomplex structures on homogeneous
manifolds},
J. reine angew. Math. {\bf 516} (1999), 159--181.
\end{thebibliography}
\end{document}